\documentclass[10pt]{amsart}
\newlength{\basicwidth}
\newlength{\shortbasicwidth}
\newlength{\basicheight}
\numberwithin{equation}{section}

\begin{document}

\title
{\bf{On a combinatorial identity of Chaundy and Bullard}}
\maketitle

\vspace{1cm}
\begin{center}
HORST ALZER$^a$ \quad{and} \quad OMRAN KOUBA$^b$
\end{center}

\vspace{1cm}
\begin{center}
$^a$ Morsbacher Stra\ss e 10, 51545 Waldbr\"ol, Germany\\
\emph{Email:} \tt{h.alzer@gmx.de}
\end{center}

\vspace{0.5cm}
\begin{center}
$^b$   Department of Mathematics, Higher Institute for Applied Sciences and Technology,\\
Damascus, Syria\\
\emph{Email:}   \tt{omran$\_$kouba@hiast.edu.sy}
\end{center}

\vspace{2.5cm}
{\bf{Abstract.}} We give two new proofs of the Chaundy-Bullard formula
$$
(1-x)^{n+1} \sum_{k=0}^m {n+k\choose k} x^k +x^{m+1}\sum_{k=0}^n {m+k\choose k} (1-x)^k=1
$$
and we prove the  ``twin formula"
$$
\frac{ (1-x)^{(n+1)}}{(n+1)!} \sum_{k=0}^m \frac{n+1}{n+k+1} \frac{ x^{(k)}}{k!}
+
\frac{ x^{(m+1)}}{(m+1)!} \sum_{k=0}^n \frac{m+1}{m+k+1} \frac{ (1-x)^{(k)}}{k!}=1,
$$
where $z^{(n)}$ denotes the rising factorial.
Moreover, we present identities involving the incomplete beta function and a certain combinatorial sum.

{\vspace{0.8cm}

\vspace{0.1cm}

\newpage

\section{Introduction}

In this paper, we study the identity
\begin{equation}
(1-x)^{n+1} \sum_{k=0}^m {n+k\choose k} x^k +x^{m+1}\sum_{k=0}^n {m+k\choose k} (1-x)^k=1,
\end{equation}
where $n$ and $m$ are nonnegative integers. Formula (1.1) was published by T.W. Chaundy and J.E. Bullard \cite{CB} in 1960. In view of its elegance, (1.1) attracted the attention of several mathematicians, who published various new proofs and generalizations. For example, the authors used generating functions, the beta integral, a hypergeometric differential equation, methods from probability theory and properties of lattics paths to establish (1.1). Moreover, we can find an $n$-variable generalization  in the literature.

In two very interesting papers on the Chaundy-Bullard identity, Koornwinder and Schlosser \cite{KS1, KS2} discussed and surveyed the different proofs and provided important information on the history of (1.1). Among others, they pointed out that Chaundy and Bullard were not the first who discovered (1.1). It was already given in implicit form by P.R. de Montmort in 1713.

The aim of this paper is to present two new proofs of the Chaundy-Bullard identity. First, we study the equation
\begin{equation}
x^{m+1} P(x)+(1-x)^{n+1} Q(x)=1,
\end{equation}
where $P$ and $Q$ are polynomials of $\deg P \leq n$ and $\deg Q \leq m$. We determine all solutions of (1.2) and apply our result
 to prove (1.1). Second, we offer a combinatorial proof of (1.1). We write  (1.1) in the form
$$
\sum_{k=0}^{n+m} A_k x^k +\sum_{k=0}^{n+m} B_k x^k=1
$$
and use certain identities for combinatorial sums to show that
$$
A_0+B_0=1, \quad A_k+B_k=0 \quad (k=1,...,m+n).
$$

These two approaches are given in Sections 2 and 3 and in Sections 4 and 5, respectively. Finally, in Section 6, we apply our results and present identities involving the incomplete beta function, the rising factorial, and we provide an identity for an alternating combinatorial 
sum with two binomial coefficients.

\vspace{0.3cm}
\section{A functional equation}

The following solution of a functional equation might be of independent interest.

\vspace{0.3cm}
{\bf{Theorem 2.1.}} \emph{Let $n,m$ be nonnegative integers.
There exist uniquely determined polynomials $P$ and $Q$ with $\deg P\leq n$, $\deg Q \leq m$ such that}
\begin{equation}
x^{m+1} P(x) +(1-x)^{n+1} Q(x)=1.
\end{equation}
\emph{These polynomials are given by}
\begin{equation}
P(x)=\sum_{k=0}^n {n+m+1\choose k} x^{n-k} (1-x)^k = \sum_{k=0}^n (-1)^k \frac{n+m+1}{k+m+1} {n+m\choose m}{n\choose k} x^k,
\end{equation}
\begin{equation}
Q(x)=\sum_{k=0}^m{n+m+1\choose k} x^k (1-x)^{m-k}=\sum_{k=0}^m {n+k\choose k} x^k.
\end{equation}

\vspace{0.3cm}
\begin{proof}
(i)  
First, we prove the uniqueness. We assume that there exist polynomials $P,Q$ and $P^*, Q^*$ with $\deg P\leq n$, $\deg P^* \leq n$ and $\deg Q\leq m$, $\deg Q^*\leq m$ such that
$$
x^{m+1} P(x)+(1-x)^{n+1} Q(x)=1
$$
and
$$
x^{m+1} P^*(x)+(1-x)^{n+1} Q^*(x)=1.
$$
Then, we obtain
$$
x^{m+1} \bigl( P(x)-P^*(x)\bigr)=(1-x)^{n+1}\bigl(  Q^*(x)-Q(x) \bigr).
$$
It follows that $x^{m+1}$ divides $(1-x)^{n+1}\bigl( Q^*(x)-Q(x)\bigr)$. Since $\mbox{gcd}(x^{m+1}, (1-x)^{n+1})=1$, we conclude that $x^{m+1}$ divides $Q^*(x)-Q(x)$ with $\deg\bigl( Q^*-Q\bigr)\leq m$. Therefore, $Q=Q^*$ and, consequently, $P=P^*$.

(ii) Let
\begin{equation}
P(x)=\sum_{k=0}^n {n+m+1\choose k} x^{n-k} (1-x)^k \quad\mbox{and} \quad Q(x)=\sum_{k=0}^m{n+m+1\choose k} x^k (1-x)^{m-k}.
\end{equation}
Then,
$$
x^{m+1} P(x)=\sum_{k=0}^n{n+m+1\choose k} x^{n+m+1-k} (1-x)^k
= \sum_{k=m+1}^{n+m+1} {n+m+1\choose k} x^k (1-x)^{n+m+1-k}
$$
and
$$
(1-x)^{n+1} Q(x)=\sum_{k=0}^m {n+m+1\choose k} x^k (1-x)^{n+m+1-k}.
$$
It follows that
$$
x^{m+1}P(x)+(1-x)^{n+1} Q(x)=\sum_{k=0}^{n+m+1} {n+m+1\choose k} x^k (1-x)^{n+m+1-k}=(1-x+x)^{n+m+1}=1.
$$

(iii) It remains to show that for the polynomials given in (2.4) we have the representations
\begin{equation}
P(x)=\sum_{k=0}^n p_k x^k
\end{equation}
with
\begin{equation}
p_k= (-1)^k \frac{n+m+1}{k+m+1} {n+m\choose m} {n\choose k} \quad(0\leq k\leq n)
\end{equation}
and
\begin{equation}
Q(x)=\sum_{k=0}^m q_k x^k
\end{equation}
with 
\begin{equation}
q_k={n+k\choose k}\quad (0\leq k\leq m).
\end{equation}
We differentiate both sides of (2.1) with respect to $x$. Then,
\begin{equation}
x^m \bigl(  (m+1)P(x)+x P'(x) \bigr) = (1-x)^n \bigl(  (n+1)Q(x)-(1-x)Q'(x) \bigr).
\end{equation}
It follows that $x^m$ divides the expression on the right-hand side of (2.9).  Since $\mbox{gcd}\bigl(  x^m, (1-x)^n \bigr)=1$, we conclude that $x^m$ divides $(n+1)Q(x)-(1-x) Q'(x)$, where
$$
\deg\bigl( (n+1) Q(x)-(1-x) Q'(x) \bigr)\leq m.
$$
It follows that there exists a constant $\mu_{n,m}$ such that
\begin{equation}
(n+1) Q(x)-(1-x) Q'(x) = \mu_{n,m} x^m.
\end{equation}
Using (2.9) we obtain
\begin{equation}
(m+1)P(x)+x P'(x)=\mu_{n,m} (1-x)^n.
\end{equation}
From (2.10) with $x=1$ we get $\mu_{n,m}=(n+1)Q(1)$ and (2.4) with $x=1$ yields $Q(1)={n+m+1\choose m}$. Thus,
\begin{equation}
\mu_{n,m}=(n+1){n+m+1\choose m}.
\end{equation}
Applying (2.7), (2.10) and (2.12)  leads to
$$
\sum_{k=0}^{m-1} \bigl(  (n+k+1)q_k -(k+1)      q_{k+1} \bigr) x^k+(n+m+1) q_m x^m=(n+1){n+m+1\choose m}x^m.
$$
Comparing the coefficients gives
$$
q_m={n+m\choose m}\quad\mbox{and}  \quad q_{k+1}=\frac{n+k+1}{k+1}q_k \quad (0\leq k\leq m-1).
$$
This yields (2.8).

Similarly, using (2.5), (2.11) and (2.12) we obtain
$$
\sum_{k=0}^n (m+k+1) p_k x^k =(n+1) {n+m+1\choose m}\sum_{k=0}^n (-1)^k {n\choose k} x^k.
$$
Comparing the coefficients leads to (2.6). The proof of Theorem 2.1 is complete.
\end{proof}

\vspace{0.3cm}
{\bf{Remark 2.2.}} The existence part of Theorem 2.1 follows from Bezout's lemma for polynomials.

\vspace{0.3cm}
\section{First proof of identity (1.1)}

We set
$$
P_{n,m}(x)=P(x) \quad\mbox{and} \quad Q_{n,m}(x)= Q(x),
$$
where $P$ and $Q$ are defined in (2.2) and (2.3), respectively. Since
$$
P_{n,m}(x)=Q_{m,n}(1-x),
$$
we obtain from (2.1),
\begin{eqnarray}\nonumber
1-(1-x)^{n+1}\sum_{k=0}^m {n+k\choose k} x^k & = & 1-(1-x)^{n+1} Q_{n,m}(x) \\ \nonumber
& = & x^{m+1} P_ {n,m}(x) \\ \nonumber
& = & x^{m+1} Q_{m,n}(1-x) \\ \nonumber
& = & x^{m+1} \sum_{k=0}^n {m+k\choose k} (1-x)^k.\nonumber
\end{eqnarray}
This leads to (1.1).

\vspace{0.3cm}
\section{Two lemmas}

In order to present our second proof of (1.1) we need the combinatorial identities given in the following two lemmas.

\vspace{0.3cm}
{\bf{Lemma 4.1.}} \emph{Let $p$ be a nonnegative integer and let $x$ be a real number. Then,}
\begin{equation}
\sum_{\nu=0}^p (-1)^{\nu} {x+\nu\choose \nu+1}{x\choose p-\nu}={x\choose p+1}.
\end{equation}

\vspace{0.3cm}
This identity is due to Brill \cite{B}, who showed that (4.1) has an application in algebraic geometry. The next lemma might be known in the literature. Since we cannot give a reference, we offer a short proof.

\vspace{0.3cm}
{\bf{Lemma 4.2.}} \emph{Let $k,n$ and $m$ be nonnegative integers with $k\leq n$. Then,}
$$
\sum_{\nu=k}^n {m+\nu\choose \nu}{\nu\choose k} =(-1)^{m}\sum_{\nu=0}^m (-1)^{\nu} {n+\nu\choose \nu}{n+1\choose m+k+1-\nu}=\frac{n+m+1}{m+k+1}{n+m\choose m}{n\choose k}.
$$

\vspace{0.3cm}
\begin{proof}
We define
$$
S_n=S_n(k,m)=\sum_{\nu=k}^n {m+\nu\choose \nu}{\nu\choose k},
\quad 
T_n=T_n(k,m)=\frac{n+m+1}{k+m+1}{n+m\choose m}{n\choose k}.
$$
Then,
$$
S_k=T_k={m+k\choose m}
$$
and
$$
S_{n+1}-S_n={n+m+1\choose m}{n+1\choose k} =T_{n+1}-T_n.
$$
By induction on $n$ we conclude that $S_n=T_n$ for all $n\geq k$.

We define
$$
R_m=R_m(k,n)=(-1)^{m}\sum_{\nu=0}^m (-1)^{\nu} {n+\nu\choose \nu}{n+1\choose m+k+1-\nu}
$$
and for nonnegative integers $\nu$,
$$
W_{\nu}=W_{\nu}(k,n,m)=(-1)^{m-\nu}\frac{\nu(\nu+n-k-m)}{(n+1)(k+m+1)}{n+
1\choose m+k+1-\nu}{n+\nu\choose n}.
$$
Then,
$$
-W_{m+1}=T_n \quad\mbox{and} \quad W_{\nu}-W_{\nu +1}=(-1)^{m+\nu}{n+\nu\choose \nu}{n+1\choose m+k+1-\nu}.
$$
It follows that
$$
-W_{m+1}=\sum_{\nu=0}^m(W_{\nu}-W_{\nu+1})=R_m.
$$
Thus, $R_m=T_n$.
\end{proof}

\vspace{0.3cm}
\section{Second proof of identity (1.1)}

We define
$$
U_{n,m}(x)=x^{m+1} \sum_{k=0}^n {m+k\choose k} (1-x)^k
\quad\mbox{and}\quad
V_{n,m}(x)=U_{m,n}(1-x)=(1-x)^{n+1} \sum_{k=0}^m {n+k\choose k} x^k.
$$
Both functions are polynomials of degree $n+m+1$. Using the binomial theorem leads to
\begin{equation}
U_{n,m}(x)=x^{m+1} \sum_{k=0}^n{m+k\choose k}\sum_{\nu=0}^k (-1)^{\nu} {k\choose \nu}x^{\nu}
=\sum_{k=m+1}^{n+m+1}a_{k-m-1} x^k
\end{equation}
with
\begin{equation}
a_k=(-1)^k \sum_{\nu=k}^n {m+\nu\choose \nu}{\nu\choose k}.
\end{equation}
Let
$$
b_k=(-1)^k {n+1\choose k} \quad(k\geq 0), \quad
c_k={n+k\choose k} \quad (0\leq k\leq m), \quad c_k=0 \quad (k\geq m+1).
$$
Applying the binomial theorem and the Cauchy product formula gives
$$
V_{n,m}(x)=\sum_{k=0}^\infty b_k x^k \sum_{k=0}^\infty c_k x^k =\sum_{k=0}^{n+m+1} d_k x^k
$$
with
$$
d_k=\sum_{\nu=0}^k b_{k-\nu} c_{\nu}.
$$
We have $d_0=1$ and for $k=1,...,m$ we get
$$
d_k
=(-1)^k \sum_{\nu=0}^k (-1)^{\nu} {n+\nu\choose \nu}{n+1\choose k-\nu}.
$$
Next, we apply  (4.1) with $p=k-1$ and $x=n+1$. We obtain
$$
0   =  {n+1\choose k}-\sum_{\nu=0}^{k-1} (-1)^{\nu} {n+1+\nu\choose \nu+1}{n+1\choose k-1-\nu} =
  \sum_{\nu=0}^{k} (-1)^{\nu} {n+\nu\choose \nu}{n+1\choose k-\nu}.
$$
This gives $d_k=0$ for $k=1,...,m$. Thus,
\begin{equation}
V_{n,m}(x)=1+\sum_{k=m+1}^{n+m+1} d_k x^k
\end{equation}
with
\begin{equation}
d_k=(-1)^k \sum_{\nu=0}^m(-1)^{\nu} {n+\nu\choose \nu}{n+1\choose k-\nu}.
\end{equation}
From (5.1) and (5.3) we conclude that
\begin{equation}
U_{n,m}(x)+V_{n,m}(x)=1+\sum_{k=m+1}^{n+m+1} (a_{k-m-1}+d_k) x^k.
\end{equation}
We set $k=l+m+1$ with $0 \leq l\leq n$. From (5.2) and (5.4) we obtain
$$
(-1)^l ( a_{k-m-1}+d_k)
=\sum_{\nu=l}^n {m+\nu\choose \nu}{\nu\choose l}-(-1)^m \sum_{\nu=0}^m (-1)^{\nu} {n+\nu\choose \nu}{n+1\choose l+m+1-\nu}.
$$
Using Lemma 4.2 gives $a_{k-m-1}+d_k=0$ for $k=m+1,...,n+m+1$, so that (5.5) yields that (1.1) is valid.

\vspace{0.3cm}
\section{Three remarks and a twin theorem}

\vspace{0.3cm}
{\bf{Remark 6.1.}} If we multiply 
 both sides of (1.1) by $x^{\alpha-1} (1-x)^{\beta-1}$, where 
$\alpha$ and $\beta$ are positive real numbers, and integrate from $x=0$ to $x=a$, then we get 
\begin{equation}
\sum_{k=0}^m {n+k\choose k} B_a(\alpha+k,\beta+n+1)+\sum_{k=0}^n {m+k\choose k} B_a(\alpha+m+1, \beta +k)=B_a(\alpha,\beta),
\end{equation}
where
$$
B_a(x,y)=\int_0^a t^{x-1} (1-t)^{y-1} dt
$$
denotes the incomplete beta function.

\vspace{0.3cm}
{\bf{Remark 6.2.}} From (2.2) and (2.3) we obtain
\begin{equation}
\sum_{k=0}^m{n+m+1\choose k}x^k (1-x)^{m-k}=\sum_{k=0}^m {n+k\choose k} x^k=(n+m+1){n+m\choose m}\sum_{k=0}^m\frac{(-1)^k}{k+n+1} {m\choose k} (1-x)^k.
\end{equation}
We multiply by $x^{\alpha-1} (1-x)^{\beta-1}$ ($\alpha,  \beta>0$) and integrate over $[0,1]$. This yields
$$
\sum_{k=0}^m{n+m+1\choose k} \frac{\Gamma(\alpha+k) \Gamma(\beta+m-k)}{\Gamma(\alpha+\beta+m)}
=\sum_{k=0}^m  {n+k\choose k}\frac{ \Gamma(\alpha+k) \Gamma(\beta)}{\Gamma(\alpha+\beta+k)}
$$
$$
= (n+m+1){n+m\choose m}\sum_{k=0}^m \frac{(-1)^k}{k+n+1}{m\choose k} \frac{\Gamma(\alpha)\Gamma(\beta+k)}{\Gamma(\alpha+\beta+k)}.
$$
Using the rising factorial notation
$$
z^{(n)}=\prod_{j=0}^{n-1} (z+j)=\frac{\Gamma(z+n)}{\Gamma(z)}
$$
 leads to
$$
\sum_{k=0}^m{n+m+1\choose k} \alpha^{(k)} \beta^{(m-k)}
=\sum_{k=0}^m  {n+k\choose k}\alpha^{(k)} (\alpha+\beta+k)^{(m-k)}
$$
$$
= (n+m+1){n+m\choose m}\sum_{k=0}^m \frac{(-1)^k}{k+n+1}{m\choose k} \beta^{(k)} (\alpha+\beta+k)^{(m-k)}
$$
which is valid for all complex numbers $\alpha$ and $\beta$ as it is a polynomial identity.

\vspace{0.3cm}
{\bf{Remark 6.3.}} Applying  the binomial theorem gives
$$
\sum_{k=0}^m\frac{(-1)^k}{k+n+1} {m\choose k} (1-x)^k=
\sum_{k=0}^m\frac{(-1)^k}{k+n+1} {m\choose k} \sum_{\nu=0}^k (-1)^{\nu} {k\choose \nu} x^{\nu}
= \sum_{k=0}^m\sum_{\nu=k}^m \frac{(-1)^{k+\nu}}{\nu+n+1} {m\choose \nu}{\nu\choose k} x^k.
$$
Comparing the coefficients we conclude from (6.2) that the identity
$$
\sum_{\nu=k}^m\frac{(-1)^{\nu}}{\nu+n+1} {m\choose \nu}{\nu\choose k}
=\frac{(-1)^k}{{n+m+1}}\frac{{n+k\choose k}}{{n+m\choose m}}
$$
is valid for all nonnegative integers $k,m,n$ with $k\leq m$.

\vspace{0.3cm}
We conclude this paper with a striking companion to (1.1), where the canonical basis $(x^k)_{k\geq 0}$ of the space of polynomials is replaced  by another basis, namely $(x^{(k)}/k!)_{k\geq 0}$.

\vspace{0.3cm}
{\bf{Theorem 6.4.}} \emph{Let $n$ and $m$ be nonnegative integers. Then,}
$$
\frac{ (1-x)^{(n+1)}}{(n+1)!} \sum_{k=0}^m \frac{n+1}{n+k+1} \frac{ x^{(k)}}{k!}
+
\frac{ x^{(m+1)}}{(m+1)!} \sum_{k=0}^n \frac{m+1}{m+k+1} \frac{ (1-x)^{(k)}}{k!}=1.
$$

\vspace{0.3cm}
\begin{proof}
Let $\alpha, \beta >0$.
Using (6.1) with $a=1$ and
$$
B(\alpha+p,\beta+q)=\frac{\Gamma(\alpha+p)\Gamma(\beta+q)}{\Gamma(\alpha+\beta+p+q)}
=\frac{\alpha^{(p)}\beta^{(q)}}{(\alpha+\beta)^{(p+q)}} B(\alpha,\beta),
$$
where $p$ and $q$ are nonnegative integers, we obtain
$$
\beta^{(n+1)} \sum_{k=0}^m{n+k\choose k} \frac{  \alpha^{(k)}}{(\alpha+\beta)^{(n+k+1)}}
+
\alpha^{(m+1)} \sum_{k=0}^n{m+k\choose k} \frac{  \beta^{(k)}}{(\alpha+\beta)^{(m+k+1)}}
=1.
$$
We set $\beta=1-\alpha$ with $\alpha\in (0,1)$. This gives
\begin{equation}
\frac{ (1-\alpha)^{(n+1)}}{(n+1)!} \sum_{k=0}^m \frac{n+1}{n+k+1} \frac{ {\alpha}^{(k)}}{k!}
+
\frac{ {\alpha} ^{(m+1)}}{(m+1)!} \sum_{k=0}^n \frac{m+1}{m+k+1} \frac{ (1-\alpha)^{(k)}}{k!}=1.
\end{equation}
By analytic continuation, we conclude that (6.3) is valid as an identity for polynomials.
\end{proof}

\section*{Conflict of Interest Declaration}

The authors declare that there is no conflict of interest.

\end{document}